\documentclass{article}

\usepackage{amssymb}
\usepackage{amsmath}
\usepackage{amsthm}

\newtheorem{theorem}{Theorem} \newtheorem{lemma}{Lemma}[section]
\newtheorem{propo}{Proposition}[section]

\def\Fo{F$\mbox{\o}$lner}
  \newcommand{\ep}{\varepsilon}
  
 \newcommand{\e}{\ep} 
 
\newcommand{\R}{\mathbb{R}} \newcommand{\Z}{\mathbb{Z}}
 \newcommand{\C}{\mathbb{C}}
 
\newcommand{\cA}{\mathcal{A}}
\newcommand{\cM}{\mathcal{M}}
\newcommand{\cN}{\mathcal{N}}
\newcommand{\cR}{\mathcal{R}} 
\newcommand{\cB}{\mathcal{B}}
\def\proof{\smallskip\noindent{\it Proof.} } 
\def\Mat{\mbox{Mat}}
\def\G{\Gamma}
\def\rk{\mbox{rk}}
\def\rkm{\mbox{rkm}}
\def\supp{\mbox{supp}}
\def\rank{\mbox{rank}}
\def\ker{\mbox{ker}}

\def\ng{\mathcal{N}(\Gamma)}\def\ug{\mathcal{U}(\Gamma)}

\title{Connes Embeddings and von Neumann Regular Closures of 
Group Algebras\footnote{
\, Research sponsored by OTKA Grant No. 69062}}

\author{G\'abor Elek}
\begin{document}
\maketitle

\begin{abstract} The analytic von Neumann regular closure $R(\Gamma)$ of a
  complex group algebra $\C\Gamma$ was introduced by Linnell and Schick.
This ring is the smallest $*$-regular subring in the algebra of affiliated
operators $U(\Gamma)$ containing $\C\Gamma$. We prove that all the algebraic
von Neumann regular closures corresponding to sofic representations of an
amenable group are isomorphic to $R(\Gamma)$. This result can be viewed as
a structural generalization of L\"uck's Approximation Theorem.

\noindent
The main tool of the proof which might be of independent interest is that
an amenable group algebra $K\Gamma$ over any field $K$ can be embedded to the
rank completion of an ultramatricial algebra.
 
\end{abstract}

\section{Introduction}

\subsection{ Regular rank rings} 
In this paper all rings are considered unital.
Regular rings were introduced by John
von Neumann, these are the rings where any principal right ideal is generated
by an idempotent (see \cite{goodearl}). 
A $*$-regular ring $\cR$ is a ring with involution and
$a^*a=0$ implies that $a=0$. In a $*$-regular ring any princial right ideal is
generated by a unique projection \cite{kaplansky}. 
A $*$-regular ring $\cR$ is {\it proper} if
$\sum^n_{i=1} a_ia_i^*=0$ implies that all of the $a_i$'s equal to $0$.
Note that $\Mat_{d\times d}(\cR)$ is regular if and only if $\cR$ is regular,
nevertheless for a $*$-regular ring $\cR$ $\Mat_{d\times d}(\cR)$ is
$\star$-regular if and only if $\cR$ is proper \cite{berberian1}. Since
$\Mat_{k\times k}(\Mat_{d\times d}(\cR))=\Mat_{kd\times kd}(\cR)$, 
$\cR$ is proper if and
only if all the matrix rings over $\cR$ are proper.

\noindent
A rank function on a regular ring $\cR$ is function 
$\rk:\cR\to\R$ satisfying the
following conditions.
\begin{enumerate}
\item $0\leq \rk(a) \leq 1$.
\item $\rk(a)=0$ if and only if $a=0$.
\item $\rk(a+b)\leq\rk(a)+\rk(b)$.
\item $\rk(ab)\leq \rk(a),\rk(b)$
\item If $e,f$ are orthogonal idempotents then $\rk(e+f)=\rk(e)+\rk(f)$.
\end{enumerate}
The most important examples of regular rank rings are matrix rings over
division rings. In this case, the values of the rank are always rational.
The rank defines a metric on the regular ring by $d(x,y)=\rk(x-y)$. The
completion of this metric is a regular rank ring as well.
 Note that for the completion
of ultramatricial algebras (see Section \ref{ultra}) the values of
 the rank can be any real number in
between zero and one \cite{goodearl}. Let $\cN$ be a finite von Neumann
algebra then its Ore localization with respect to its non-zero divisiors
 $U(\cN)$ is a $\star$-regular ring. The
elements of this ring are called affiliated operators (see \cite{reich}).
 The rank of an
affiliated operator is the trace of the idempotent that
generates the right ideal generated by the operator. Note that if $A\in
U(\cN)$,then
\begin{equation}
\label{rankformula}
\rk(A)=1-\lim_{t\to\infty}\int_0^t tr_N (E_\lambda)d\lambda\,,
\end{equation}
where $\int^\infty_0 E_\lambda d\lambda$ is the spectral decomposition
of the unbounded operator $A^*A$. This shows that if $i:\cN\to\cM$ is
a trace-preserving homomorphism between finite von Neumann algebras, then
its Ore-extension $\tilde{i}:U(\cN)\to U(\cM)$ is a rank preserving
$\star$-homomorphism.
\noindent
Note that $U(\cN)$ is
always proper (see Section \ref{linsection}). 

\noindent
If $\cR$ is regular rank ring then there is a unique natural extension of the
rank to a matrix rank of $\Mat_{k\times k}(\cR)$ \cite{halperin}.
 Note that a matrix rank
$\rkm$ has the same property as the rank $\rk$ except that 
$0\leq \rkm(M)\leq k$.

\subsection{The Connes Embedding Problem}

Let $\nu=\{d_1< d_2 <\dots\}$ be an infinite sequence of positive integers.
Then one can consider the ultraproduct of the matrix algebras
$\{\Mat_{d_i\times d_i}(\C)\}^\infty_{i=1}$ as tracial
algebras the following way (see \cite{pestov}).

Let $\omega$ be a nonprincipal ultrafilter on the natural numbers and let
$\lim_\omega$ be the corresponding ultralimit.
First, consider the algebra of bounded elements
$$\cB=\{(a_1,a_2,\dots)\in \prod_{i=1}^\infty \Mat_{d_i\times d_i}(\C) 
 \,\mid\, \sup \|a_i\|<\infty\}\,.$$
Now let $\mathcal{I}\lhd\cB$ be the ideal of elements $\{a_i\}^\infty_{i=1}$
such that $\lim_\omega\frac{tr(a_n^*a_n)}{d_n}=0$. Then
$\cB/\mathcal{I}=\cM_\nu$ is a type $II_1$-von Neumann factor with trace
defined the following way.
$$Tr_\omega[\{a_i\}^\infty_{i=1}]=\lim_\omega\frac{tr(a_n)}{d_n}\,.$$

The following conjecture is generally referred to as the Connes Embedding
Problem. Is it true that a type-$II_1$-von Neumann algebra with a separable 
predual
have a trace-preserving embedding to some $\cM_\nu$ ? See the survey of
Pestov \cite{pestov} for further details.

\noindent
There is a purely algebraic version of the Connes Embedding Problem first
considered in \cite{elekszabo}. Namely, we can consider the ultraproduct of
the matrix rings $\{\Mat_{d_i\times d_i}(\C)\}$ as {\it rank algebras}.
Let $\mathcal{J}\lhd \prod_{i=1}^{\infty} \Mat_{d_i\times d_i}(\C)$ be the
following ideal,
$$\mathcal{J}=\{\,\{a_i\}^\infty_{i=1} \,\mid\, lim_\omega\frac
{rank(a_n)}{d_n}=0\}\,.$$
Then $\prod_{i=1}^{\infty} \Mat_{d_i\times d_i}(\C)/\mathcal{J}=\cM_\nu^{alg}$
 is a simple
complete $\star$-regular rank ring. One can ask of course, whether any 
countable
dimensional regular rank ring embeds to some $\cM_\nu^{alg}$.
\subsection{L\"uck's Approximation Theorem}
Let $\Gamma$ be a finitely generated residually finite group and
let \\ $\Gamma=N_0\supset N_1\supset N_2\dots, \cap^\infty_{k=1}N_k=\{1\}$
be finite index normal subgroups. Let $\Delta\in \Mat_{d\times d}(\Z\Gamma)$
be a $d\times d$-matrix over the integer group algebra $\Z\Gamma$.
Denote by $\ng$ the von Neumann algebra of $\Gamma$. Note that $\Delta$ acts
on $(l^2(\Gamma)^d)$ as a bounded operator. Then one can
define $\dim_\Gamma Ker(\Delta)$, the von Neumann dimension of the kernel of
$\Delta$. 
Let $\pi_k:\C\Gamma\to \C(\Gamma\backslash N_k)$ the natural projection.
That is $\pi_k(\Delta)\in \Mat_{d\times d}(\C(\Gamma\backslash N_k))$ is a
finite dimensional linear transformation. According to L\"uck's Approximation
Theorem (see \cite{lueck}) \begin{equation} \label{appro}
\lim_{k\to\infty}\frac{\dim_\C Ker (\pi_k(\Delta))}{|\Gamma:N_k|}=\dim_\G
  Ker (\Delta) \end{equation}
It is conjectured that (\ref{appro}) holds for any $\Delta\in \Mat_{d\times
  d}(\C\Gamma)$ as well. The conjecture was confirmed for amenable groups
$\Gamma$ in \cite{elekappro}.
\subsection{Regular Closures}
Linnell and Schick \cite{linnellschick} proved the following theorem 
(see Section \ref{linsection}).
Let $\cR$ be a proper $\star$-regular ring. Then for any subset $T\subseteq
\cR$ there exists a smallest $\star$-regular subring containing $T$. We call
this ring $R(T,\cR)$ the regular closure of $T$ in $\cR$.
Let $\Gamma$ be a countable group, then one can consider the natural embedding
of its complex group algebra to its von Neumann algebra
$\C(\Gamma)\to\cN(\Gamma)$. Let $\ug$ be Ore localization of $\cN(\Gamma)$.
Then  $\ug$ is a proper $\star$-regular ring (see \cite{berberian2}. 
Therefore one can
consider the {\it analytic regular closure} $R(\C(\Gamma),\ug)=R(\Gamma)$.

\noindent
Now let $\Gamma=N_0\rhd N_1\rhd\dots, \cap^\infty_{i=1} N_i=\{1\}$ be normal
subgroups of a residually finite group.Let $\pi_i:\C\Gamma\to \C(\Gamma/N_i)$
be the natural projection as in the previous subsection and let
$s_i:\C(\Gamma/N_i)\to \Mat_{\Gamma/N_i\times \Gamma/N_i}(\C)$ be the natural
  representation by convolutions.
Define $r_i=s_i\circ\pi_i:\C\Gamma\to
  \Mat_{\Gamma/N_i\times \Gamma/N_i}(\C)\,.$ Then we have an injective (see
    \cite{elekszabo}) $\star$-homomorphism $r:\C\Gamma\to \cM^{alg}_\nu$, where
$\nu=\{|\Gamma/N_1|, |\Gamma/N_2|,\dots\}$.
Therefore we can consider the {\it algebraic regular closure}
$R(\C\Gamma,\cM^{alg}_\nu)$ for any normal chain of residually finite group.
The main result of this paper is the following theorem.
\begin{theorem} \label{main} Let $\Gamma$ be a finitely generated amenable
  group. Then there is a rank preserving $\star$-homomorphism 
$$j:R(\Gamma)\to R(\C\Gamma,\cM^{alg}_\nu)$$\
which is the identity map restricted on $\C\Gamma$. \end{theorem}
This theorem can be viewed as a structural generalization of  L\"uck's
Approximation Theorem for amenable groups. 
Indeed, let $\Delta\in \Mat_{k\times k}(\C\Gamma)$ then the
approximation theorem is equivalent to the fact that
$$\rkm_1(\Delta)=\rkm_2(\Delta)\,,$$
where $\rkm_1$ resp. $\rkm_2$ are the matrix rank on $\Mat_{k\times
  k}(U(\Gamma))$
resp. on $\Mat_{k\times k}(\cM^{alg}_\nu)$. However, both $\rkm_1(\Delta)$ and
$\rkm_2(\Delta)$ are equal to the matrix rank of $\Delta$ in
$\Mat_{k\times k} (R(\Gamma))$. Actually, we prove a generalization of 
Theorem \ref{main}, where we consider algebraic closures associated to
arbitrary sofic representations (see Section \ref{sofic}) of 
the group $\Gamma$. 
We shall also prove the following theorem. \begin{theorem}
\label{tetel2}
For any finitely generated amenable group and coefficient field $K$,
 the group algebra $K\Gamma$ embeds to the rank completion of an
 ultramatricial algebra. \end{theorem}
\section{von Neumann regular closures}
\label{linsection}
In this section, we review some results of Linnell and Schick \cite{linnell},
\cite{linnellschick} about the von Neumann regular closures in proper
$\star$-regular rings.
The starting points of Linnell's paper are the following two observations
about finite von Neumann algebras already mentioned in the introduction.
\begin{enumerate}
\item $U(\cN)$ is a
$\star$-regular ring, that is any right ideal is generated by a single
projection.
\item If $\alpha,\beta\in U(\cN)$ and $\alpha\alpha^*+\beta\beta^*=0$ then
$\alpha=\beta=0\,.$
\end{enumerate}
Although Linnell and Schick consider only group von Neumann algebras all what
they used are the two properties above.
The following result is a strengthening of the second observation.
\begin{propo}[Lemma 2. \cite{linnell}]
If $\alpha,\beta\in U(\cN)$ then $(\alpha\alpha^*+\beta\beta^*)U(\cN)\supseteq
\alpha U(\cN)$.
\end{propo}
Using a simple induction one also has the following proposition.
\begin{propo} [Lemma 2.5 \cite{linnellschick}]
If $\alpha_1, \alpha_2,\dots, \alpha_n\in U(\cN)$ then
$$\sum^n_{i=1} \alpha_i\alpha^*_i U(\cN)\supseteq \alpha_1 U(\cN)\,.$$
\end{propo}
This leads to the crucial proposition about the existence of the von Neumann
regular closures.
\begin{propo}[Proposition 3.1 \cite{linnellschick}] \label{proplins}
Let $\{R_i\,\mid\, i\in I\}$ be a collection of $\star$-regular
subrings of $U(\cN)$. Then $\cap_{i\in I} R_i$ is also a $\star$-regular
subring of $U(\cN)$.
\end{propo}
We also need to show that the proposition above holds for $\cM^{alg}_\mu$
as well. 
\begin{propo}
 Let $\{R_i\,\mid\, i\in I\}$ be a collection of $\star$-regular
subrings of $\cM^{alg}_\mu$. Then $\cap_{i\in I} R_i$ is also a $\star$-regular
subring of $\cM^{alg}_\mu$.
\end{propo}
\proof
Since $\cM^{alg}_\mu$ is a $\star$-regular ring we only need to prove
that if $\alpha\alpha^*+\beta\beta^*=0$ in $\cM^{alg}_\mu$, then both
$\alpha$ and $\beta$ equal to $0$. Then the proof of Proposition \ref{proplins}
works without any change.
\begin{lemma}
For finite dimensional matrices $A,B\in Mat_{k\times k}(\C)$
$$\rank(AA^*+BB^*)\geq \max(\rank (A), \rank( B))$$
\end{lemma}
\proof If $(AA^*+BB^*)(v)=0$ then $A^*(v)=0$ and $B^*(v)=0$. Hence
$$\ker (AA^*+BB^*)\subseteq \ker(A^*)\cap \ker(B^*)\,.$$
Therefore
$\rank(AA^*+BB^*)\geq\rank(A^*)=\rank(A)$ and $\rank(AA^*+BB^*)\geq\rank(B^*)
=\rank(B)\,$ \quad\qed

\vskip0.2in
Now let $A_n, B_n\in \Mat_{d_n\times d_n}(\C)$, then
$$\lim_{\omega} \frac{\rank(A_nA^*_n+B_nB_n^*)}{d_n}\geq
\lim_{\omega} \frac{\rank(A_n)}{d_n}$$
and
$$\lim_{\omega} \frac{\rank(A_nA^*_n+B_nB_n^*)}{d_n}\geq
\lim_{\omega} \frac{\rank(B_n)}{d_n}\,.$$
Hence the proposition follows.

\section{Bratteli Diagrams, Ultramatricial Algebras and Tilings}
\label{ultra}
Recall that a Bratteli diagram is an oriented countable graph such that the
vertex set is partitioned into finite sets $\{Z_i\}^\infty_{i=1}$ such a way
that

\begin{itemize}
\item If the starting vertex of an edge is $Z_i$, then the end vertex is
  necessarily in $Z_{i+1}$.
\item Each vertex has at least one outgoing edge.
\item Each vertex $\alpha$ has a non-negative {\it size} $S(\alpha)$.
\item Each edge (from a vertex $\alpha$ to a vertex $\beta$) has a
  non-negative multiplicity $K(\alpha,\beta)$ such that for each
$\beta\in Z_{n+1}$, $S(\beta)=\sum_{\alpha\in Z_n} S(\alpha) K(\alpha,\beta)$
\end{itemize}

Let $P_n$ be a probability 
distribution function on $Z_n$. We call the system $\{P_n\}^\infty_{n=1}$ a
{\it harmonic function} if
$$P_n(\alpha)=\sum_{\beta\in Z_{n+1}}\frac{ S(\alpha) K(\alpha,\beta)}
{S(\beta)} P_n(\beta)$$
for any $n\geq 1$ and $\alpha\in Z_{n}$.
\noindent
An {\it ultramatricial algebra} is constructed the following way. For each
$n\geq 1$ one consider a product ring $\oplus_{l=1}^{i_n}\Mat_{d^n_l\times
  d^n_l}(\C)$. Let $K(d^n_l,d^{n+1}_j)$ be non-negative integers satisfying
$$d^{n+1}_j=\sum_{1\leq l\leq i_n} d^n_j K(d^n_l,d^{n+1}_j)$$
for any $n\geq 1$ and $1\leq j\leq i_{n+1}$.

Now for any $n\geq 1$ and $1\leq l\leq i_{n+1}$
choose a diagonal embedding
$$E_{n,l}:\oplus^{i_n}_{l=1}(\Mat_{d^{n}_l\times
  d^n_l}(\C))^{K(d^n_l,d^{n+1}_j)}\to \Mat_{d^{n+1}_j\times d^{n+1}_j}(\C)\,.$$
The embeddings define injective maps 
$$\phi_n:\oplus^{i_n}_{l=1} \Mat_{d^{n}_l\times
  d^n_l}(\C)\to 
\oplus^{i_{n+1}}_{j=1} \Mat_{d^{n+1}_j\times
  d^{n+1}_j}(\C)\,.$$
The direct limit $\lim_{\to}\phi_n$ is the ultramatricial algebra
$\cA_\phi$. Clearly, $\cA_\phi$ is a $\star$-regular ring.

\noindent
Now for any $n\geq 1$ and $1\leq l \leq i_n$ let $P(d^n_l)$ be real numbers 
satisfying
\begin{equation}
\sum^{i_n}_{l=1} P(d^n_l)=1
\end{equation}
and
\begin{equation} P(d^n_l)=\sum_{j=1}^{i_{n+1}}
\frac {d^n_l K(d^n_l,d^{n+1}_j)}{d^{n+1}_j}
P(d^{n+1}_j)\,.
\end{equation}
Then we have a Bratteli diagram with a harmonic function, where
the vertices in $Z_n$ are exactly $\{\Mat_{d^n_l\times d^n_l}(\C)\}^{i_n}_{l=1}
$, with sizes
$\{d^n_l\}^{i_n}_{l=1}$. 
The Bratteli diagram defines a rank function $\rk_\phi$ on $\cA_\phi$.
Namely, let 
$$\rk_\phi(a_1\oplus a_2\oplus\dots\oplus a_{i_n})=
\sum^{i_n}_{l=1} m(a_l) \frac{\rank (a_l)}{d^n_l},$$
where $m(a_l)=P(d^n_l)$ and $\rank(a_l)$ is the rank of the matrix $a_l$. 
Then it is easy to see that each
$\phi_n$ is a rank preserving $\star$-isomorphism. Therefore $\cA_\phi$ is a
rank regular ring.

\noindent
Now let $\Gamma$ be a finitely generated group with a symmetric generating
system $S$. The Cayley graph of $\Gamma$, $Cay(\Gamma,S)$ is defined as
follows.  
\begin{itemize}
\item $V(Cay(\Gamma,S))=\Gamma$
\item $(a,b)\in E(Cay(\Gamma,S))$ if $as=b$ for some $s\in S$.
\end{itemize}
Let $F\subset\Gamma$ be a finite set. Then $\partial F$ is the set of vertices
that are adjacent to a vertex in the complement of $F$. The isoperimetric
constant of $F$ is defined as
$$i(F):=\frac{|\partial F|}{|F|}\,.$$
The group $\Gamma$ is amenable if there exists a \Fo-sequence in $\Gamma$
that is a a sequence of finite sets $\{F_n\}^\infty_{n=1}$ such that
$\lim_{n\to\infty} i(F_n)=0\,.$

\noindent
Now we define Bratteli-tiling systems. If $\gamma\in\Gamma,F\subset\Gamma$
then $\gamma F$ is called a $F$-tile.
A Bratteli system has the following properties.
\begin{itemize}
\item The level set $Z_n$ consists of finite sets $F^n_1,
  F^n_2,\dots,F_{i_n}^n$  and the set $E_n$ containing only the unit element.
Also, we have $i(F^n_j)\leq \frac{1}{2^n}$ for all $j$ and $n$.
\item For any $n\geq 2$ and $F^n_j\in Z_n$ we have a partition
  $F^n_j=\cup_{i=1} ^{a_{n,j}} \gamma_iA_i$, where $A_i\in Z_{n-1}$. That is
we have tiling of $F^n_j$ with the tiles of $Z_{n-1}$.
\item $K(F^{n-1}_l,F^n_j)$ is the number of $F^{n-1}_l$-tiles in the partition
  of $F^n_j$. Also $K(E_{n-1},F^n_j)$ is the number of $E_{n-1}$-tiles (single
  vertices).
\item $S(F^n_j)=|F^n_j|$, $S(E_n)=1$.
\item We also suppose that $K(E_{n-1},F^n_j)\leq\frac{1}{2^{n-1}}|F^n_j|\,.$
\end{itemize}
Let $m:\cup_{n=1}\infty Z_n\to\R$ be a harmonic function such that
$m(E_n)\to 0$ as $n\to\infty$. Then we call a system above a Bratteli tiling
system. Our main technical tool is the following proposition.
\begin{propo}\label{bratt1}
For any amenable group $\Gamma$ and generating system $S$ we can construct
a Bratteli tiling system with the following property.
For any $\epsilon>0$ and $n>0$ there exist $\delta=\delta_{\e,n}>0$ such that
if $F\in\Gamma$ is a finite set and $i(F)<\delta$ then one can tile $F$ with 
translates of the elements $Z_n$ satisfying the following property.
If $L\in Z_n$ and $T^F_L$ is the set of points in $F$ covered by a translate
of $L$ then
$$\left|\frac{|T^F_L|}{|F|}-m(L)\right|<\epsilon\,.$$
\end{propo}
\section{Proof of Proposition \ref{bratt1}}
First, let us recall the notion of $\e$-quasitilings. Let $Cay(\Gamma,S)$ be
the Cayley-graph of an amenable group $\Gamma$ as above. Let $F\subset\Gamma$
be a finite set and $A_1,A_2,\dots, A_n$ be subsets of $F$. We say that
$\{A_i\}^n_{i=1}$ $\e$-cover $F$ if
$$\frac{|\cup^n_{i=1} A_i|}{|F|}>1-\e\,.$$
Also, we call $\{A_i\}^n_{i=1}$ $\e$-disjoint if there exist disjoint sets
$\{B_i\}^n_{i=1}$, $B_i\subset A_i$, such that
$$\frac{|B_i|}{|A_i|}>1-\e\,.$$
The system $\{A_i\}^n_{i=1}$ $\e$-quasi-tiles $F$ if it both $\e$-covers $F$ and
$\e$-disjoint.
The following result of Ornstein and Weiss \cite{ornsteinweiss} is crucial for our proof.
\begin{propo} [Quasitiling theorem]\label{quasi}
Let $F_1\subset F_2\subset\dots$ be a \Fo-sequence. Then for any $\e>0$
ther exists $\delta>0$ and a subfamily $F_{n_1}\subset F_{n_2} \subset\dots
F_{n_k}$ such that if $i(F)<\delta$ then $F$ can be $\e$-quasitiled by
translates of the $F_{n_i}$'s.\end{propo}
Observe that if $i(A)<\e$ and $B\subset A$, $\frac{|B|}{|A|}>1-\e$, then
$$i(B)<(d+1)\frac{\e}{1-\e}\,,$$
where $d$ is the degree of the vertices of $Cay(\Gamma,S)$. Indeed, $\partial
B$
is covered by the union of $\partial A$ and the neighbours of the vertices in
$A\backslash B$. Thus 
$|\partial B|\leq |\partial A| + 
d |A\backslash B|$. Therefore,
$$\frac{|\partial B|}{|B|}\leq \frac{|\partial B|}{|A|(1-\e)}
\leq (d+1) \frac{\e}{1-\e}\,.$$
Now using the quasitiling theorem we construct a Bratteli system inductively.
Suppose that $\{F^m_1, F^m_2,\dots, F^m_{i_m}\}$ and the decreasing 
sequence of positive constants 
$\{\delta_n\}^m_{i=0}$ are 
already given such a way
that
\begin{itemize}
\item for any $i\geq 1$ $|\partial F^m_i|<\min(\delta_{n_1},\frac{1}{2^n})$,
\item if $i(F)<\delta_m$ then $F$ can be tiled by translates of the $F^m_i$'s
and less than $(1/2^m)|F|$ single points.
\end{itemize}
Now let $G_1\subset G_2\subset\dots$ be a \Fo-sequence and $c>0$ such that
if $B\subset G_j$ for some $j$ and $\frac{|B|}{|G_j|}>1-c$ then
$$i(B)<\min (\delta_n,1/2^{n+1})\,.$$
By Proposition \ref{quasi} there exists a family of finite subsets
$F^{n+1}_1, F^{n+1}_2,\dots, F^{n+1}_{i_{n+1}}$ (namely subsets of a certain
system $G_{n_1}, G_{n_2},\dots, G_{n_k}$) and a constant $\delta_{n+1}$
 such that
\begin{itemize}
\item for any $i\geq 1$ $|\partial F^{n+1}_i|
<\min(\delta_n,\frac{1}{2^{n+1}})$,
\item if $i(F)<\delta_{n+1}$ then $F$ can be tiled by translates of the 
$F^{n+1}_i$'s
and less than $(1/2^{n+1})|F|$ single points.
\end{itemize}

By the induction above one can obtain a Bratteli system. Now we construct a
harmonic function $m$.
Fix a \Fo-sequence $H_1\subset H_2\subset\dots$. Let
$\{\delta_n\}^\infty_{n=1}$ be the constans as above.
If
$$\min(\frac{1}{j_{n+1}},\delta_{j_{n+1}})\leq i(H_n)< 
\min(\frac{1}{j_n},\delta_{j_n})$$
then pick a tiling of $H_n$ by translates of
the elements of $Z_{j_n}$ such a way that the number of
single vertices is less than $(\frac{1}{2^{j_n}})|H_n|$.
Then pick a tiling of the 
$Z_{j_n}$-tiles by translates of of $Z_{j_n-1}$ such a way that the number of
single vertices in any tile $T$ is less than $(\frac{1}{2^{j_{n-1}}})|T|$.
Inductively, we obtain a tiling of $H_n$ by $Z_i$-translates for any $1\leq i
\leq j_n$.
Note that the number of single vertices used in the $Z_i$-tiling of $H_n$
is less than 
$$(\sum_{k=i}^{j_n} \frac{1}{2^k})|H_n|\leq \frac{1}{2^{i-1}}|H_n|\,.$$
 If $A\in Z_i$ then denote by $c_k(A)$ the number of vertices in
$H_k$ covered by $A$-translates and let
$$m_k(A)=\frac{c_k(A)}{|H_k|}\,.$$
Clearly, $\sum_{A\in Z_i} m_k(A)=1\,.$ We may suppose that
for any $A$, $\lim_{k\to\infty} m_k(A)=m(A)$ exists, otherwise we could pick
a subsequence of $\{H_k\}^\infty_{k=1}$.
\begin{lemma}
The function $m$ is harmonic satisfying $\lim_{i\to\infty} m(E_i)=0\,.$
\end{lemma}
\proof
The fact that $\lim_{i\to\infty} m(E_i)=0\,.$ follows from our previous
observation about the number of single vertices used in the tiling
of $H_k$. By definition, if $A\in Z_i$
$$m_{k}(A)=\sum_{B\in Z_{i+1}}\frac{S(A) K(A,B)}{S(B)} m_{k}(B)\,.$$
By taking the limit as $n\to\infty$ we get that

$$m (A)=\sum_{B\in Z_{i+1}}\frac{S(A) K(A,B)}{S(B)} m (B)\,.\qed $$
Now let us show that our Bratteli tiling system satisfies the required
property. First we prove a simple lemma.
\begin{lemma} \label{l42}
For any $i>0$ and $\delta>0$ there exists $\lambda>0$ and $p>0$ with 
the following
property. Let $k>p$ and $J\subseteq H_{k}$, $\frac{|J|}{|H_{k}|}>1-\lambda\,.$
For $A\in Z_i, |A|>1$, let $j_A^{k}$ be the number of vertices
in $H_{k}$ that are covered by an $A$-translate (in the tiling previously 
defined) which is completely inside $J$. Also, let $j^{k}_{E_i}$
be the number of points in $H_{k}$ that are not covered by any of the
$A$-translates above. Then 
$$\left|\frac{j_A^{k}}{H_k}-m(A)\right|<\delta\quad\mbox{and}\quad 
\left|\frac{j^{k}_{E_i}}{H_k}-m(E_i)\right|<\delta\,.$$
\end{lemma}
\proof The number of points covered by such $A$-translates that contain at
least one point from the complement of $J$ is less than $|H_{k}\backslash J|
|A|\,.$ Hence
$$m_{k}(A)\geq j_A^{k}\geq m_{k}(A)-\frac{|H_{k}\backslash J|}
{|H_{k}|}|A|\,.$$
Since $\sup_{A\in Z_i}|A|<\infty$ and $m_{k}(A)\to m(A)$ the lemma easily 
follows. \qed

\vskip 0.2in
\noindent
Now let $\e$ be the constant in our proposition and
$0<\alpha <\e/2$. By Proposition \ref{quasi}, we have a subfamily
$H_{a_1},H_{a_2},\dots,H_{a_t}$ of $\{H_{k}\}^\infty_{k=1}$ that
$\alpha$-quasitiles any finite set $F$ with $i(F)<\delta_\alpha$.
By the previous lemma it
 means that we have disjoint subsets $J$ in $F$ that can be tiled
by $Z_i$-translates such a way that using the notation of our proposition
$$\left|\frac{T^F_A}{|F|}-m(A)\right|<\frac{\e}{10}\,,$$
for any $A\in Z_i$ provided that $\alpha$ is small enough. 
We cover all the remaining points (that are not in the
$J$'s) by single vertices. Then we get the required tiling of $F$. \qed

\section{The canonical rank on amenable group algebras}
In this section we recall some results from \cite{elekamena}.
Let $\Gamma$ be a finitely generated amenable group and $K\Gamma$ be its group
algebra over the field $K$. Let $\{F_n\}^\infty_{n=1}$ be a \Fo-
sequence. 
For $a\in K\Gamma$ let $V^a_n\subset K^{F_n}\subset K\Gamma$ be the vector
space of elements $z$ supported on $F_n$ such that $za=0$. Then
$$\lim_{n\to\infty}\frac{\dim_K V^a_n}{|F_n|}=k_a$$
exists and independent on the choice of the \Fo-sequence.
We call $\rk(a)=1-k_a$ the natural rank of $a$. It is proved in 
\cite{eleklinnell}
that if $K=\C$ then $\rk(a)=1-\dim_\Gamma Ker\,M_a$, where $dim_\Gamma$ is the
von Neumann dimension and $Ker\,M_a$ is the set of elements 
$w\in l^2(\Gamma)$ for which
$wa=0$. Note that the rank can be computed slightly differently as well.
Let $S$ be a symmetric generating system of $\Gamma$ and $Cay(\Gamma,S)$ be
the Cayley-graph of $\Gamma$.
We consider the shortest path metric $d_{Cay(\Gamma,S)}$ on $\Gamma$.
Let $\supp(a)\subset B_r(1)$, where $B_r(1)$ is the $r$-ball around the unit
element in the Cayley-graph and
$$\supp(a)=\{\gamma\in\Gamma\,\mid\,a_\gamma\neq 0\}$$
if $a=\sum a_\gamma \gamma$. For a finite set $F\subset \Gamma$, let
$\partial_r F$
be the set of elements $x\in F$ such that
$$d_{Cay(\Gamma,S)}(x,F^c)\leq r\,.$$
Note that $\partial F=\partial_1(F)$.
Clearly, if $b\in K\Gamma$ and
$\supp(b)\in F\backslash \partial_r F$ then $\supp (ba)\subset F$.
Then for any $s>r$
\begin{equation}\label{rankf}
\rk(a)=\lim_{n\to\infty} \frac{\dim_K W_na} {|F_n|},
\end{equation}
where $W_n$ is the set of elements in $K\Gamma$
supported on $F_n\backslash \partial_s F_n$.
\section{Sofic representations}\label{sofic}
\subsection{Sofic approximation and sofic representations} \label{soficsub1}
In this section we recall the notion of sofic representations from 
\cite{elekszabo}.
Let $\Gamma$ be a finitely generated group with a symmetric generating set
$S$. Let $\{G_n\}^\infty_{n=1}$ be a sequence of finite graphs such that
the directed edges are labeled by the elements of $S$ such a way that if
$(x,y)$ is labeled by $s$, then $(y,x)$ is labeled by $s^{-1}$.
We say that $\{G_n\}^\infty_{n=1}$ is a sofic approximation of $\Gamma$ if for
any natural number $r>0$ there exists $n_r>0$ such that 
\begin{itemize}
\item if $n\geq n_r$ then for the set $V^r_n$ of vertices $x$ for which
the ball $B_r(x)$ in $G_n$ is isomorphic to the ball $B_r(1)\in Cay(\Gamma,S)$
as labeled graphs
$$ \frac{|V^r_n|}{|V(G_n)|}>1-\frac{1}{r}\,.$$
\end{itemize}
A group is called {\it sofic} if it possesses a sofic approximation. In this 
moment no non-sofic group is known. If $Cay(\Gamma,S)$ is the Cayley-graph on
an amenable group and $\{F_n\}^\infty_{n=1}$ is a \Fo-sequence then the
induced graphs $G_n$ of the sets $F_n$ form a sofic approximation of $\Gamma$.
If $\Gamma$ is residually finite (amenable or not) with normal chain 
$\{N_k\}^\infty_{k=1},\,\cap^\infty_{k=1} N_k=\{1\}\,$ then the
graph sequence $Cay(\Gamma/N_k,S)$ form a sofic approximation of $\Gamma$.
If $\{G_n\}^\infty_{n=1}$ is an arbitrary sofic approximation of a group
$\Gamma$ then one can construct an imbedding of $K\Gamma$ ($K$ is an arbitrary
field) to the ultraproduct of matrix algebras the following way.

\noindent
Let $\{\Mat_{V(G_n)\times V(G_n)}(K)\}^\infty_{n=1}$ be a sequence of
matrix algebras. Let $a\in K\Gamma$, $a=\sum r_\gamma \gamma$ be an element
of the group algebra such that if $r_\gamma\neq 0$ then $\gamma\in
B_r(1)\subset Cay(\Gamma,S)$.
Let $\{e_x\}_{x\in V(G_n)}$ be the natural basis of $K^{V(G_n)}$. If $x\in
V^r_n$ then let 
$$\psi_n(a)(e_x)=\sum_{y\in B_r(x)} k_ye_y\,,$$
where $k_y=r_\gamma$ if $x\gamma=y$. Note that by our condition on the
support of $a$ $x\gamma =y$ is meaningful. If $x\notin V^r_n$ let
$\psi_n(a)(x)=0$.
This way one can define an injective homomorphism
$\psi_\mu:K\Gamma\to\cM^{alg}_\mu$, where $\mu=\{|V(G_1)|, |V(G_2)|,\dots\}$
 If $K=\C$ the homomorphism above is a
$*$-homomorphism. The map $\psi_\mu$ is called the sofic representation
associated to the sequence $\{G_n\}^\infty_{n=1}$. 
\subsection{Sofic approximation of amenable groups}
Let $\{G_n\}^\infty_{n=1}$ be a sofic approximation of the amenable
group $\Gamma$ (with symmetric generating set $S$). For $L>0$ let $Q^{G_n}_L$
be the set of vertices $x$ in $G_n$ such that 
$$B_L(x)\cong B_L(1)\subset Cay(\Gamma, S)$$
as $S$-labeled graphs.
If $F\subset B_L(1)$ then for $x\in Q^{G_n}_L$ we call $\pi(F)$ an
$F$-translate, where $\pi:B_L(1)\to B_L(x)$ is the $S$-labeled graph
isomorphism mapping $1$ to $x$.
In \cite{elekappro} we proved the following generalization of the
Ornstein-Weiss quasitiling theorem.
\begin{propo} \label{soficweiss}
Let $F_1\subset F_2\dots$ be a \Fo-sequence. Then for any $\e>0$ there
exists $L>0$, $\delta>0$ and a finite subcollection $F=\{F_{n_1},
F_{n_2},\dots,F_{n_t}\}$ such that if
$$\frac{Q^{G_n}_L}{|V(G_n)|}>1-\delta $$
then $G_n$ can be $\e$-quasitiled by $F$-translates.
\end{propo}
\section{Imbedding $K\Gamma$ to the completion of an ultramatricial algebra}
Let $(\{Z_i\}^\infty_{i=1},m)$ be the Bratteli tiling system as in Proposition
\ref{bratt1}.
We construct an ultramatricial algebra as in Section \ref{ultra}.
Let $\oplus_{A\in Z_i} \Mat_{|A|\times |A|}(K)$ be the $i$-th algebra. 
For $B\in Z_{i+1}$ let
$$E_B:\oplus_{A\in Z_i} \Mat_{|A|\times |A|}(K)\to \Mat_{|B|\times |B|}(K)$$
be the diagonal embedding, where the image of each $\Mat_{|A|\times |A|}(K)$
is $K(A,B)$ $|A|\times |A|$-diagonal block in $ \Mat_{|B|\times |B|}(K)$.
Let 
$$\phi_i=\oplus_{B\in Z_{i+1}}E_B:\Mat_{|A|\times |A|}(K)\to 
\oplus_{B\in Z_{i+1}}\Mat_{|B|\times |B|}(K)$$
the product map.
Now we define the maps
$\pi^A_i:K\Gamma\to \Mat_{|A|\times |A|}(K)$ the following way.
We identify the elements of $\Mat_{|A|\times |A|}(K)$ with the linear
transformations from $K^A$ to $K^A$ the natural way.
Let $a\in K\Gamma$, $\supp (a)\subset B_r(1)$.
If $x\in A\backslash \partial_r(A)$, then let
$$\pi^A_i(a)(e_x)=\sum a_\gamma e_{x\gamma}\,.$$
Note that by the condition on the support $x\gamma$ is well-defined.
If $\partial_r(A)$, then let $\pi^A_i(a)(e_x)=0\,.$
Finally we define the maps $\pi_i:=\oplus_{A\in Z_i} 
\pi^A_i:K\Gamma\to \oplus_{A\in Z_i} 
\Mat_{|A|\times |A|}(K)$.
\begin{lemma}
\label{cauchy}
For any $a\in K\Gamma$, $\{[\pi_i(a)]\}^\infty_{i=1}$ is a Cauchy-sequence in
$\cA_\phi$, where $[\pi_i(a)]$ denotes the image of $\pi_i(a)$ under the natural\\
embedding $\oplus_{A\in Z_i} \Mat_{|A|\times |A|}(K)\to \cA_\phi$.
\end{lemma}
\proof
First of all note that
$$\rk_\phi(\phi_i\circ \pi_i(a)-\pi_{i+1}(a))=
\sum_{B\in Z_{i+1}}m(B) \frac{\rank(E_B\circ\pi_i(a)-\pi_{i+1}^B(a))}{|B|}\,.$$

Observe that
$$\rank(E_B\circ\pi_i(a)-\pi^B_{i+1}(a))=
|B|-\dim_K \ker(E_B\circ\pi_i(a)-\pi^B_{i+1}(a)).$$
On the other hand,
$$\dim_K \ker (E_B\circ\pi_i(a)-\pi^B_{i+1}(a))\leq T_B\,,$$
where $T_B$ is the number of vertices in $B$ for which
$$E_B\circ\pi_i(a)(e_x)=\pi^B_{i+1}(a)(e_x)\,.$$
Clearly,
$$T_B\geq |B|-|\partial_r B|-\sum_{A\in Z_i} K(A,B) |\partial_r(A)|\,.$$
Recall that if $|B|>1$ then
$|\partial B|\leq |B| 2^{-(i+1)}\,.$
Hence $$|\partial_r B|\leq |B|2^{-(i+1)} (d+1)^{r+1}\,,$$
where $d$ is the degree of the vertices in the Cayley graph. 
Also, 
$$\sum_{A\in Z_i} K(A,B) |\partial_r(A)|\leq
K(E_i,B)+\sum_{A\in Z_i\,,|A|>1} K(A,B)|A| 2^{-i}(d+1)^{r+1}\,.$$
Therefore,
$$T_B\geq |B|-|B|2^{-(i+1)}(d+1)^{r+1}-|B| 2^{-(i+1)}-|B|2^{-i}(d+1)^{r+1}\,.$$
Hence,
$$\rk_\phi(\phi_i\circ \pi_i(a)-\pi_{i+1}(a))\leq $$$$\leq
2^{-(i+1)} +\sum_{B\in Z_{i+1}\,,|B|>1}m(B) (2^{-(i+1)}(d+1)^{r+1} +
2^{-(i+1)} + 2^{-i}(d+1)^{r+1})$$
Thus the lemma follows.\qed
\begin{lemma}
Let $a\in K\Gamma, b\in K\Gamma$, then
\begin{enumerate}
\item $lim_{i\to\infty} \rk_\phi(\phi_i(a)\phi_i(b)-\phi_i(ab))=0\,.$
\item $lim_{i\to\infty} \rk_\phi(\phi_i(a)+\phi_i(b)-\phi_i(a+b))=0\,.$
\item If $K=\C$ then
$lim_{i\to\infty} \rk_\phi(\phi_i(a^*)-\phi_i^*(a))=0\,.$
\end{enumerate}\end{lemma}
\proof We prove the first part only, the other two statements can be seen
exactly the same way.
If $x\in A\backslash \partial_{r+s}(A)$ then 
$$\phi_i(a)\phi(b)-\phi_i(ab)(e_x)=0\,.$$
Therefore 
$$\rk_\phi(\phi_i(a)\phi_i(b)-\phi_i(ab))
\leq \frac{|\partial_{r+s}(A)|}{|A|}\,.$$
Hence the lemma follows. \qed

\vskip 0.2in
\noindent
Let $\phi(a)\in\overline{\cA}_\phi$ be the limit of elements $\phi_i(a)$.
By the previous lemma $\phi$ is a homomorphism and if $K=\C$ then $\phi$ is
even a $\star$-homomorphism. Finally, we prove that
$\rk_\phi(\phi(a))=\rk_\Gamma(a)\,.$
By definition,
$$\rk_\phi(\phi_i(a))=\sum_{A\in Z_i}m(A)\frac{\rank(\phi_i(a))}{|A|}\,.$$
If $|A|=1$, then $m(A)\leq 2^{-i}$ otherwise by (\ref{rankf})
$$\lim_{i\to\infty}\frac{\rank(\phi_i(a))}{|A|}=\rk_\Gamma(a)\,.$$\
Hence, $\rk_\phi(\phi(a))=\rk_\Gamma(a)\,.$ This finishes the proof of our
theorem. \qed
\section{The proof of the main theorem}
\subsection{The strategy of the proof}
We have four complete regular $*$-rings: $\overline{\cA_\phi}$, $U(\Gamma)$,
$\cM^{alg}_\mu$ and $U(\cM_\mu)$. Also, we define seven rank preserving
embeddings
\begin{enumerate}
\item $f_1:\C\Gamma\to \cM^{alg}_\mu$
\item $f_2:\C\Gamma\to U(\Gamma)$
\item $f_3:\C\Gamma\to U(\cM_\mu)$
\item $f_4:\C\Gamma\to \overline{\cA_\phi}$
\item $f_5:\overline{\cA_\phi}\to\cM^{alg}_\mu$
\item $f_6:\overline{\cA_\phi}\to U(\cM_\mu)$
\item $f_7:U(\Gamma)\to U(\cM_\mu)$.
\end{enumerate}
We shall show three identities:
\begin{enumerate}
\item $f_5\circ f_4=f_1$
\item $f_6\circ f_4=f_3$.
\item $f_7\circ f_2= f_3$
\end{enumerate}
From these identities the main theorem easily follows. Indeed,
$R(\Gamma)$ is the smallest $\star$-regular ring containing $\C\Gamma$ in
$U(\Gamma)$. The ring $R(\Gamma)$ is inside $\overline{\cA_\phi}$, in fact,
it is the minimal $\star$-regular ring containing  $\C\Gamma$
in $\overline{\cA_\phi}$. On the other hand,
$R(\C\Gamma, \cM^{alg}_\mu)$ is also the smallest $\star$-regular ring
 containing  $\C\Gamma$ in $\overline{\cA_\phi}$. \qed
\subsection{The first identity}\label{sub81}
Let $\{G_n\}^\infty_{n=1}$ be the sofic approximation of our group $\Gamma$
and
$\cM^{alg}_\mu$ be the associated ultraproduct. Let $\{H_n\}^\infty_{n=1}$
be the \Fo-sequence in the proof of Proposition \ref{bratt1}.
Let $f_4$ be the map $\phi:\C\Gamma\to \overline{\cA_\phi}$ defined in
the proof of Theorem \ref{tetel2}. Let $f_1$ be the map $\psi_\mu:
\C\Gamma\to \cM^{alg}_\mu$ defined in Subsection \ref{soficsub1}.
Fix $k\geq 1$. Now we define maps
$\tau_{k,n}:\oplus_{A\in Z_k} \Mat_{A\times A}(\C)\to \Mat_{V(G_n)\times
  V(G_n)}(\C)$
for large enough $n\geq 1$. 

\noindent
First, let $q\geq 1$ be an integer. We say that $G_n$ is $q$-regular if
$G$ can be $\frac{1}{2^q}$-quasitiled by translates of
$\{H_{n_1}, H_{n_2},\dots H_{n_t}\}$, where $n_i\geq q$. For $n\geq 1$, let
$q(n)$ be the largest $q$ for which $G_n$ is $q$-regular. By Proposition
\ref{soficweiss}, for any $q\geq 1$ there exists some $n_q$ such that
if $n\geq n_q$ then $q(n)\geq q$.

\noindent
Now consider a $\frac{1}{2^q}$-quasitiling of $G_n$ by translates of
$\{H_{n_1}, H_{n_2},\dots H_{n_t}\}$. Then consider the iterated tiling for
each $H_{n_i}$ above by $Z_k'$s as in the proof of Proposition \ref{bratt1},
starting with $Z_{l(n)}$-tilings. Since the translates are not disjoint this 
does not yet define a tiling of $G_n$. However, let $\{J_\alpha\}_{\alpha\in
  I}$ be the disjoint parts in the quasitiling. That is each $J_\alpha$ is
inside some $H_{n_i}$-translate having size at least $(1-\frac{1}{2^q})
|H_{n_i}|$. Discard the tiles that are inside some $Z_{l(n)}$-tile that
is not contained in some $J_\alpha$. Cover, the remaining part of $G_n$ by
single vertices.
For $A\in Z_k$, let $Q_n(A)$ be the number of vertices in  $V(G_n)$ 
that are covered
by an $A$-translate. By Lemma \ref{l42}, it is easy to see that
$$\lim_{n\to\infty}\frac{Q_n(A)}{|V(G_n)|}=m(A)\,.$$
Now let $\tau_{k,n}:\oplus_{A\in Z_k} \Mat_{A\times A}(\C)\to 
\Mat_{V(G_n)\times V(G_n)}(\C)$ be the natural diagonal map induced by the
tiling above. If $|A|>1$, the definition is clear. The case where  $A=E_k$, 
that is $A$ is a
single point needs some clarification. In the diagonal map, we use only
those vertices in $G_n$ that are in some ``good'' $Z_{l(n)}$-tile, in other
words, that are not used to cover the remaining part. All the diagonal
elements in the image of $\tau_{k,n}$ that belong to a vertex covering the
remaining part are zero.

\noindent
By the iterative tiling construction, one can immediately see that
$\tau_{k,n}\circ \phi_k=\tau_{k+1,n}$. If $k>q(n)$, let us define
$\tau_{k,n}:=0$. Therefore we have a map
$$\tau=(\tau_1,\tau_2\dots):\oplus_{A\in Z_k}\Mat_{A\times A}(\C)\to
\oplus^\infty_{n=1} \Mat_{V(G_n)\times V(G_n)}(\C)$$
and this map extends to $\cA_\phi$ as well. 
\begin{lemma}\label{ranklemma}
For any $(a_1\oplus a_2\oplus\dots a_{i_k})\in \oplus_{A\in Z_i}
\Mat_{A\times A}(\C)$
$$\lim_{n\to\infty} \frac{\rank\, 
(\tau_{k,n}(a_1\oplus a_2\oplus\dots a_{i_k}))}
{|V(G_n)|}=\rk_\phi(a_1\oplus a_2\oplus\dots a_{i_k})\,.$$
\end{lemma}
\proof
$$\lim_{n\to\infty} \frac{\rank (\tau_{k,n}(a_1\oplus a_2\oplus\dots a_{i_k}))}
{|V(G_n)|}=\lim_{n\to\infty} \sum_{A\in Z_k}
\frac{\frac {Q_n(A)}{|A|} \rank(a_A)} {|V(G_n)|}=$$
$$=\sum_{A\in Z_k} m(A)\frac{rank(a_A)}{|A|}=
\rk_\phi(a_1\oplus a_2\oplus\dots
a_{i_k})\,.$$
\qed
\vskip 0.2in

\noindent
Therefore we have a rank-preserving map
 $\tau_{alg}:\oplus_{A\in Z_k}\Mat_{A\times A}(\C)\to \cM^{alg}_\mu$
defined as $\pi\circ\tau_{alg}$, where 
$$\pi:\oplus^\infty_{n=1} \Mat_{V(G_n)\times V(G_n)}(\C)\to \cM^{alg}_\mu$$
is the quotient map. The map $\tau_{alg}$ extends to the rank completion of
$\cA_\phi$, resulting in the map $f_5$.

\noindent
Now let us prove the first identity.
Let $a\in\C\Gamma$, $\supp(a)\subset B_r(1)\subset Cay(\Gamma,S)$.
Then $f_1(a)$ can be represented in $\oplus^\infty_{n=1} 
\Mat_{V(G_n)\times V(G_n)}(\C)$ by the element $\oplus_{n=1}^\infty \psi_n(a)$,
where $\psi_n$ is defined in Subsection  \ref{soficsub1}.
On the other hand, $f_5\circ f_4(a)$ is represented by
$\oplus_{n=1}^\infty \psi'_n(a)$, where
$$\psi'_n(a)(e_x)=\sum_{y\in B_r(x)} k_ye_y\,,$$
where $k_y=r_\gamma$ if $x\gamma=y$ and $x\in\partial_r(J_\alpha)$, for some
$J_\alpha$ in a $H_{n_i}$-translate. Clearly,
$$\lim_{n\to\infty} \frac{z_n(a)}{|V(G_n)|}\,,$$
where $z_n(a)$ is the number of elements $x\in V(G_n)$ for which
$$\psi'_n(a)(e_x)=\psi_n(a)(e_x)\,.$$
Therefore $f_5\circ f_4=f_1$.
\subsection{The second identity}

Let $\rk_1$ resp. $\rk_2$ denote the ranks on $\cM_\mu$ resp. $\cM_\mu^{alg}$
Let $$t\in\oplus^\infty_{n=1} \Mat_{V(G_n)\times V(G_n)}(\C)=(t_1,
t_2,\dots)\,,$$
where $\sup \|t_i\|<\infty$. Note that $t$ represents and element
$[t]_{\cM_\mu}$ in $\cM_\mu$ and an element
$[t]_{\cM_\mu^{alg}}$ in $ \cM_\mu^{alg}$. It is important to note
that $rk_1([t]_{\cM_\mu})$ is not necessarily equal to 
$rk_2([t]_{\cM_\mu^{alg}})$. Indeed, let $t_n=\frac{1}{n} Id$. Then
$\rk_1([t]_{\cM_\mu})=0\,.$ Nevertheless, $rk_2([t]_{\cM_\mu^{alg}})=1\,.$
However, we have the following lemma.
\begin{lemma}\label{however1}
Let $t=(t_1,t_2,\dots)\in\oplus^\infty_{n=1} \Mat_{V(G_n)\times V(G_n)}(\C)$,
where for any $n\geq 1$, $t_n$ is  self-adjoint and all the $t_n$'s have
altogether finitely many distinct eigenvalues $\lambda_0=0,
\lambda_1,\lambda_2,\dots,\lambda_k$. Then
$$\rk_1([t]_{\cM_\mu})=\rk_2([t]_{\cM^{alg}_\mu})\,.$$
\end{lemma}
\proof Let $t_{n,i}$ be the multiplicity of $\lambda_i$ in $t_n$. Then
$$\rk_2([t]_{\cM^{alg}_\mu})=\lim_\omega(1-\frac{t_{n,0}}{|V(G_n)|})\,.$$
The spectral decomposition of $[t]_{\cM_\mu}$ is $\sum^k_{i=1}\lambda_i P_i$,
where $$tr_{\cM_\mu}(P_i)=\lim_\omega (1-\frac{t_{n,i}}{|V(G_n)|})\,.$$
By (\ref{rankformula}) $$\rk_1([t]_{\cM_\mu})=
\lim_\omega(1-\frac{t_{n,0}}{|V(G_n)|})\,.\quad\qed$$
We also need the following lemma.
\begin{lemma} \label{however2}
Let $t$ be as above and suppose that
$$\lim_{n\to\infty}\frac{\rank(t_n)}{|V(G_n)|}=0\,.$$ Then
$[t]_{\cM_\mu}=0\,.$
\end{lemma}
\proof We need to check that
$\lim_{n\to\infty}\frac{tr (t_n^* t_n)}{|V(G_n)|}=0\,.$
Observe that $\sup \|t_n^*t_n\|=K<\infty$ and $\rank(t^*_nt_n)\leq
K\rank(t^*_nt_n)\,.$
Then $tr(t^*_nt_n)\leq K\rank(t^*_nt_n)$. Hence the lemma follows. \qed
\vskip 0.2in
\noindent
Let $i_\mu:\C\Gamma\to\cM_\mu$ be defined as $\rho\circ\psi$, where
$\psi=\oplus^\infty_{n=1}\psi_n$ as in Subsection \ref{soficsub1} and
$$\rho:B(\oplus^\infty_{n=1}\Mat_{V(G_n)\times V(G_n)}(\C))\to\cM_\mu\,$$
be the quotient map from the bounded part of the direct product. The map
$i_\mu$ is trace-preserving and extends to an injective trace-preserving map
$\overline{i_\mu}:\cN(\Gamma)\to\cM_\mu$ (see \cite{Elekhiper} and
\cite{pestov}). The map $f_3$ is just $i_\mu$ composed by the embedding of
$\cM_\mu$ into its Ore-extension. We prove that $f_3$ is rank-preserving later.

\noindent
Now let us define the map $f_6$. Let $\tau$ be the map 
defined in Subsection \ref{sub81}. Then let
$j:\rho\circ\tau:\cA_\phi\to\cM_\mu$ and let $s=u\circ j$, where $u:\cM_\mu\to
U(\cM_\mu)$ be the natural embedding. Then $f_6$ is defined as the extension
of $s$ onto $\overline{\cA_\phi}$. We need to show of course that $j$ is 
rank-preserving that is
$$\rk_1[\tau(\underline{a})]_{\cM_\mu}=
\rk_2[\tau(\underline{a})]_{\cM^{alg}_\mu}\,,$$
for any $\underline{a}\in\oplus_{A\in Z_k} \Mat_{A\times A}(\C)\,.$
However, this immediately follows from Lemma \ref{however1}.

\noindent
Now we prove the second indentity. This also shows that $f_3$ is
rank-preserving. Again, it is enough to show that
\begin{equation}
\label{utolso}
[\oplus^\infty_{n=1}\psi_n(a)-\psi'_n(a))]_{\cM_\mu}=0\,.
\end{equation}
We already proved that
$$\lim_{n\to\infty}\frac{\rank(\psi_n(a)-\psi'_n(a))}{|V(G_n)|}=0\,.$$
Obviously, $\sup \|\psi_n(a)-\psi'_n(a)\|<\infty\,,$
hence (\ref{utolso}) follows from Lemma \ref{however2}.
\subsection{The third identity}
By definition, $i_\mu=\overline{i_\mu}\circ i$, where $i$ is the natural
embedding of $\C\Gamma$ into $\ng$. This immediately proves the third
identity. This completes the proof of our main theorem.\qed


\begin{thebibliography}{9}

\bibitem{berberian1} {\sc S. K. Berberian}, Baer $*$-rings. 
Die Grundlehren der mathematischen Wissenschaften, Band 195. 
Springer-Verlag, New York-Berlin, (1972)

\bibitem{berberian2} {\sc S. K. Berberian}, The maximal ring of quotients
of a finite von Neumann algebra. {\sl Rocky Mount.
J. Math. } {\bf 12} (1982) no. 1, 149--164.

\bibitem{elekamena} {\sc G. Elek},
The rank of finitely generated modules over group algebras.
{\em Proc. Amer. Math. Soc.} {\bf 131} (2003) no. 11, 3477--3485.

\bibitem{Elekhiper} {\sc G. Elek and E. Szab\'o}, Hyperlinearity, essentially
 free actions and $L^2$-invariants. The sofic property. {\sl  Math. Ann.}  
{\bf 332}  (2005),  no. 2, 421--441. 

\bibitem{eleklinnell} {\sc G. Elek}, On the analytic zero 
divisor conjecture of Linnell. {\em Bull. Lond. Math. Soc} {\bf 35}
(2003) no 2., 236--238.

\bibitem{elekszabo} {\sc G. Elek and E. Szab\'o}, 
Sofic groups and direct finiteness.
{\em  Journal of Algebra} 280 (2004) no.2, 426--434.

\bibitem{elekappro} {\sc G. Elek}, The strong approximation conjecture holds 
for amenable groups. {\em Journal of Funct. Anal.} {\bf 239} (2006) no. 1,
345--355.

\bibitem{goodearl} {\sc K. R. Goodearl}, von Neumann regular rings.
{\sl Robert E. Krieger Publishing Co.}, Inc., Malabar, FL, (1991)

\bibitem{halperin} {\sc I. Halperin}, Extensions of the rank function.
{\sl Studia Math.} {\bf 27} (1966) 325--335.

\bibitem{kaplansky} {\sc I. Kaplansky}, Rings of operators.
{\sl W. A. Benjamin}, Inc., New York-Amsterdam (1968)

\bibitem{linnell} {\sc P. Linnell}, Embedding group algebras into finite
von Neumann regular rings. {\sl Modules and Comodules, Trends in Mathematics,
Birkhauser Verlag}, 295--300.

\bibitem{linnellschick} {\sc P. Linnell and T. Schick},
The Atiyah conjecture and Artinian rings
{\sl electronic} http://arxiv.org/pdf/0711.3328

\bibitem{lueck} {\sc W. L\"uck}, Approximating $L\sp 2$-invariants by their
 finite-dimensional analogues. {\sl Geom. Funct. Anal.} {\bf  4}
 (1994)  no. 4, 455--481.
\bibitem{ornsteinweiss} {\sc D. S. Ornstein} and {\sc B. Weiss}, 
Entropy and isomorphism
theorems for actions of amenable groups, {\sl J. Anal. Math} {\bf 48}
(1987) 1-141.
\bibitem{pestov} {\sc  V. G. Pestov}, Hyperlinear and sofic groups: a brief
guide. {\sl Bull. Symb. Logic} {\bf 14} (2008), no. 4, 449--480.

\bibitem{reich} {\sc H. Reich}, 
 On the $K$- and $L$-theory of the algebra of operators affiliated to a 
finite von Neumann algebra.  {\sl K-theory} {\bf 24} (2001) no.4 303--326.

\bibitem{thom} {\sc A. Thom}, $L^2$-cohomology for von Neumann algebras.
{\sl GAFA} {\bf 18} (2008) no. 1 251--270.











\end{thebibliography}
\end{document}